\documentclass[reqno,centertags,12pt]{amsart}
\usepackage{amsmath,amsthm,amscd,amssymb,latexsym,verbatim}

\usepackage{graphicx,epsf,cite}

-\marginparwidth \oddsidemargin -4mm \evensidemargin \oddsidemargin

\newtheorem{theorem}{Theorem}[section]

\newtheorem{lemma}[theorem]{Lemma}

\theoremstyle{definition}

\theoremstyle{remark}

\newcounter{smalllist}

\allowdisplaybreaks
\numberwithin{equation}{section}

\newcommand{\lb}{\label}

\newcommand{\beq}{\begin{equation}}
\newcommand{\eeq}{\end{equation}}

\newcommand{\bal}{\begin{align}}
\newcommand{\eal}{\end{align}}
\newcommand{\bals}{\begin{align*}}
\newcommand{\eals}{\end{align*}}

\newcommand{\bbR}{{\mathbb{R}}}

\newcommand{\bbZ}{{\mathbb{Z}}}

\newcommand{\bbT}{{\mathbb{T}}}

\newcommand{\eps}{\varepsilon}

\newcommand{\til}{\tilde}

\begin{document}
\title[Exponential Growth for Euler Equations]
{Exponential Growth of the Vorticity Gradient \\ for the Euler Equation on the Torus}

\author{Andrej Zlato\v s}

\address{\noindent Department of Mathematics \\ University of
Wisconsin \\ Madison, WI 53706, USA \newline Email: \tt
zlatos@math.wisc.edu}


\begin{abstract} 
We  prove that there are solutions to the Euler equation on the torus with $C^{1,\alpha}$ vorticity and smooth except at one point such that 
the vorticity gradient grows in $L^\infty$ at least exponentially as $t\to\infty$.  
The same result is  shown to hold for 
the vorticity Hessian and 
smooth solutions.  Our proofs use a version of a recent result by Kiselev and \v Sver\' ak \cite{KS}.
\end{abstract}

\maketitle

\section{Introduction} \lb{S1}

Let $(2\bbT)^2=[-1,1)^2$ be the two-dimensional torus (i.e., we identify opposite sides of the square) and consider the Euler equation on $(2\bbT)^2$, in vorticity formulation:
\beq\lb{1.1}
\omega_t+u\cdot\nabla\omega = 0, \qquad \omega(0,\cdot)=\omega_0.
\eeq
The velocity $u$ is found from the vorticity $\omega$ via the Biot-Savart law
\[
u(t,x)=  \frac 1{2\pi} \sum_{n\in\bbZ^2} \int_{[-1,1]^2}  \frac{(x_2-y_2-2n_2,-x_1+y_1+2n_1)}{|x-y-2n|^2} \omega(t,y)dy, 
\]
obtained by taking the Biot-Savart kernel $K(x)=\tfrac 1{2\pi}(x_2,-x_1)|x|^{-2}$ on $\bbR^2$ and extending $\omega$ periodically.
Initial data $\omega_0$ will here be $C^1$ and odd in both $x_1$ and $x_2$, hence the latter property will hold for all $t\ge 0$, as well as
\[
\|\omega(t,\cdot)\|_{L^\infty} = \|\omega_0\|_{L^\infty}.
\]

Global regularity of bounded solutions to \eqref{1.1} was first proved by Wolibner \cite{Wol} and H\" older \cite{Hol}.
We consider here the question of how fast the gradient of $\omega$ can grow in $L^\infty$ as $t\to\infty$.  The well known upper bound is double-exponential $\|\nabla \omega(t,\cdot)\|_{L^\infty}\le e^{Ce^{Ct}}$ but it has been a long standing open question whether this is attainable.   The best infinite time result in the plane or on the torus (i.e., domains without a boundary) so far has been the proof of the possibility of super-linear growth for smooth solutions by Denisov \cite{Den}.  He also proved that double-exponential growth is possible on arbitrarily long finite time intervals \cite{Den2}, and constructed patch solutions to the 2D Euler equation with a regular prescribed stirring for which the distance between two approaching patches decreases double-exponentially in time \cite{Den3}.  

For domains with boundaries (and no flow boundary condition),  Yudovich \cite{Yud1, Yud2} and Nadirashvili \cite{Nad} provided examples with unbounded growth and at least linear growth, respectively.
These results have been dramatically improved in a striking recent work by Kiselev and \v Sver\' ak \cite{KS}, who proved the possibility of infinite time double-exponential growth of  the vorticity gradient in a disc, thereby answering in the affirmative the above open question in this setting.   The boundary is crucial in \cite{KS} and the double-exponential growth is proved to occur on it as well.    

Related to this, we note that the double-exponential upper bound is only known for 2D domains with regular boundaries.  In fact, Kiselev and the author \cite{KZ} proved that there are domains whose boundary  is smooth except at two points, and on which some solutions to the Euler equation with smooth initial data blow up in finite time.
We refer the reader to \cite{KS} for more history and further references related to \eqref{1.1}.  

In the present paper we prove that on the torus, at least  exponential growth of the vorticity gradient happens for some  $C^{1,\alpha}$ initial data, as well as that such growth is possible for the vorticity Hessian for smooth initial data.
Our proof uses a sharper version of a key result from \cite{KS} (see Lemma \ref{L.2.1} below), applied  on the torus instead of the disc.

\begin{theorem}\lb{T.1.1}
(i) For any $\alpha\in(0,1)$ and $A<\infty$, there is $\omega_0\in C^{1,\alpha}((2\bbT)^2)$ with $\|\omega_0\|_{L^\infty}\le 1$ and there is $T_0\ge 0$ such that the  solution of \eqref{1.1} satisfies for all $T\ge T_0$,
\[
\sup_{t\le T} \|\nabla \omega(t,\cdot)\|_{L^\infty([0,2\exp(-AT)]^2)} \ge e^{AT}.
\] 

(ii)  For any $A<\infty$, there is $\omega_0\in C^{\infty}((2\bbT)^2)$ with $\|\omega_0\|_{L^\infty}\le 1$ and there is $T_0\ge 0$ such that the  solution of \eqref{1.1} satisfies for all $T\ge T_0$,
\[
\sup_{t\le T} \|D^2 \omega(t,\cdot)\|_{L^\infty([0,2\exp(-AT)]^2)} \ge e^{AT}.
\] 

\end{theorem}


Our $\omega_0$ will be very simple.  For instance, in (i) it can be smooth except at the origin, with $\omega_0(r,\phi)=r^{1+\alpha}\sin(2\phi)$ (in polar coordinates) near the origin.  Then $\omega$ remains in $C^{1,\alpha}$, and since $u(t,0)=0$ for all $t>0$ by symmetry (oddness of $\omega$ in $x_1,x_2$), it follows that both $u,\omega$ are smooth except at the origin at all times.  

We also note that we will have $\|\omega_0\|_{L^\infty}=1$, but as in \cite{KS}, this can be arbitrarily small.


\smallskip

{\bf Acknowledgements.}  
The author thanks Sergey Denisov and Alexander Kiselev for useful discussions and comments.  He also 
acknowledges partial support  by NSF grants DMS-1056327  and DMS-1159133.

\section{Proof of Theorem \ref{T.1.1}} \lb{S2}

For $x\in[0,1]^2$ we let $\til x:=(-x_1,x_2)$, $\bar x:=(x_1,-x_2)$, and $Q(x):=[x_1,1]\times[x_2,1]$.
In what follows, $C$ will always be some universal constant which may change between inequalities.

\begin{lemma}\lb{L.2.1}
Let  $\omega(t,\cdot)\in L^\infty((2\bbT)^2)$  be odd in both $x_1$ and $x_2$.  If  $x_1,x_2\in[0, \tfrac 12]$, then
\beq\lb{2.1}
u_j(t,x)=   (-1)^j \left( \frac 4\pi\int_{Q(2x)}\frac{y_1y_2}{|y|^4} \omega(t,y) dy 
+  B_j(t,x)    \right) x_j  \qquad (j=1,2)
\eeq
where, with some universal $C$,
\begin{align}
|B_1(t,x)|\le C \|\omega(t,\cdot)\|_{L^\infty} \left( 1 + \min \left\{ \log \left( 1+ \frac {x_2}{x_1} \right) , x_2 \frac{\|\nabla\omega(t,\cdot)\|_{L^\infty([0,2x_2]^2)}}{\|\omega(t,\cdot)\|_{L^\infty}} \right\} \right), \lb{2.7}
\\ |B_2(t,x)|\le C\|\omega(t,\cdot)\|_{L^\infty} \left( 1 + \min \left\{ \log \left( 1+ \frac {x_1}{x_2} \right), x_1 \frac{\|\nabla\omega(t,\cdot)\|_{L^\infty([0,2x_1]^2)}}{\|\omega(t,\cdot)\|_{L^\infty}} \right\} \right). \lb{2.8}
\end{align}  
\end{lemma}

{\it Remark.}  
If $c x_1\ge x_2$ (resp.~$cx_2\ge x_1$) for some $c<\infty$, then the $\min$ can obviously be dropped in \eqref{2.7} (resp.~\eqref{2.8}) as long as $C=C(c)$.  This version of these formulas was proved in \cite{KS} on the disc.  In that case $Q(2x)$ can be replaced by $Q(x)$ as well, which is done in \cite{KS}.

\begin{proof}
Let us only consider 
$j=1$ because $j=2$ follows by symmetry: $K$ shows that if $\til \omega(t,x):=\omega(t,x_2,x_1)$, then $\til u(t,x)=-(u_2(t,x_2,x_1),u_1(t,x_2,x_1))$.  
The Biot-Savart law gives
\beq \lb{2.5}
u_1(t,x) 
=  \frac 2{\pi} \sum_{n\in\bbZ^2} \int_{[0,1]^2}  \left[ \frac{y_1(x_1-2n_1)(x_2-y_2-2n_2)}{|x-(y+2n)|^2 |x-(\til y+2n)|^2}
-\frac{y_1(x_1-2n_1)(x_2+y_2-2n_2)}{|x-(\bar y+2n)|^2 |x-(- y+2n)|^2}  \right] \omega(t,y) dy, 
\eeq
by  using the symmetries $\omega(t,\til y)=-\omega(t,y)$ and then $\omega(t,\bar y)=-\omega(t,y)$ to express the integral over $[-1,1]^2$ via that over $[0,1]^2$.  

Let us first consider the right-hand side of \eqref{2.5} with the term $n=(0,0)$ removed.  
The  first term in the integral equals (recall that $x_1,x_2\in[0,\tfrac 12]$)
\[
\frac{-2y_1n_1(x_2-y_2-2n_2)}{|x-(y+2n)|^2 |x-(\til y+2n)|^2} + x_1 O(|n|^{-3}).
\]
We combine it with the same term for $\til n=(-n_1,n_2)$ to obtain 
\[
\frac{-32y_1n_1(x_2-y_2-2n_2)x_1[n_1(x_2-y_2-2n_2)^2+x_1^2-y_1^2+4n_1^2]}{|x-(y+2n)|^2 |x-(\til y+2n)|^2|x-(y+2\til n)|^2 |x-(\til y+2 \til n)|^2} + x_1 O(|n|^{-3}) = x_1 O(|n|^{-3}).
\]
This means that
\[
\left| \frac 2\pi \sum_{n\neq(0,0)} \int_{[0,1]^2}  \frac{y_1(x_1-2n_1)(x_2-y_2-2n_2)}{|x-(y+2n)|^2 |x-(\til y+2n)|^2}
 \omega(t,y) dy \right| \le C x_1\|\omega(t,\cdot)\|_{L^\infty}.
\]
An identical argument proves this also for the second term in  \eqref{2.5}.  

We therefore only need to consider the term with $n=(0,0)$, which is $x_1$ times
\beq\lb{2.4}
\frac 2\pi \int_{[0,1]^2} \left[ \frac{y_1(x_2-y_2)}{|x-y|^2 |x-\til y|^2}
-\frac{y_1(x_2+y_2)}{|x-\bar y|^2 |x+ y|^2}  \right] \omega(t,y) dy. 
\eeq
We will show that this equals to $\tfrac 4\pi$ times the integral in \eqref{2.1}, plus an error controlled by the right-hand side of  \eqref{2.7}, thus proving \eqref{2.1}.  We will again only consider the first term in the integral since the second will be handled in the same way.

We  separate the integral into either 3 or 4  regions.  If $x_1\ge x_2$, then these regions will be $Q(2x), [0,2x_1]\times[0,1]$, and $[2x_1,1]\times[0,2x_2]$.  If $x_1<x_2$, we also split the last region into $[2x_1,2x_2]\times[0,2x_2]$ and $[2x_2,1]\times[0,2x_2]$.  The 3 region case  is parallel to the treatment in \cite{KS} (where the domain is a disc).
In the 4 region case we 
need to obtain an extra estimate for 
the integral over $[2x_1,2x_2]\times[0,2x_2]$ 
(this is not necessary for the second term in \eqref{2.4}).

We start with $Q(2x)$, where  we have
\[
\int_{Q(2x)}  \frac{y_1x_2}{|x-y|^2 |x-\til y|^2} dy 
\le C \int_{Q(2x)}  \frac{|x|}{|y|^3} dy \le C \int_{|2x|}^1 \frac{|x|}{r^2} dr \le C.
\]
In $Q(2x)$ also
\[
\frac{y_1y_2}{|x-y|^2 |x-\til y|^2} - \frac{y_1y_2}{|y|^4 } = \frac{O(|x||y|^5)}{|y|^8} = O(|x||y|^{-3}), 
\]
so the integral of the absolute value of this difference over $Q(2x)$ is also bounded by
\[
C\int_{|2x|}^1 \frac{|x|}{r^2} dr \le C.
\]
Hence integration of the first term in \eqref{2.4} over $Q(2x)$ gives ($\tfrac 2\pi$ times) the integral in \eqref{2.1},  plus an error  bounded by $C\|\omega(t,\cdot)\|_{L^\infty}$.  (Integration of the second term in \eqref{2.4} gives the same, whence the factor of $\tfrac 4\pi$ in \eqref{2.1}.)  

Next integrate over $[0,2x_1]\times[0,1]$. After substituting $z_j:=y_j-x_j$, the absolute value of the integral of the first term in \eqref{2.4} can be bounded by $\|\omega(t,\cdot)\|_{L^\infty}$ times
\[
C\int_{0}^{x_1} \int_{0}^{1}  \frac{x_1 z_2}{(z_1^2+z_2^2) (x_1^2+z_2^2)} dz_2dz_1 
= C\int_{0}^{1}  \frac{x_1}{ x_1^2+z_2^2} \arctan \frac {x_1}{z_2} dz_2 \le C\frac\pi 2\arctan \frac {1}{x_1} \le C.
\]

Finally, the corresponding integral over $[2x_1,1]\times[0,2x_2]$ is bounded by $\|\omega(t,\cdot)\|_{L^\infty}$ times
\[
C \int_{x_1}^1 \int_{0}^{x_2}  \frac{z_1 z_2}{(z_1^2+z_2^2)^2} dz_2dz_1 
\le C \int_{0}^{x_2}  \frac{z_2}{x_1^2+z_2^2}  dz_2 
\le C  \log \frac{x_1^2+x_2^2}{x_1^2}  \le C\log\left( 1+\frac{x_2}{x_1} \right).
\]
This gives the first term in the $\min$ in \eqref{2.7}.  If $x_1\ge x_2$, then we are done because the min is $\le 1$ and can be absorbed in the 1 in \eqref{2.7}.  

So let us assume  $x_1<x_2$, and perform the above integration over $[2x_2,1]\times[0,2x_2]$ instead of $[2x_1,1]\times[0,2x_2]$.  We see that the integral is bounded by $C \|\omega(t,\cdot)\|_{L^\infty} \log (1+\tfrac{x_2}{x_2}) =C \|\omega(t,\cdot)\|_{L^\infty} \log 2$, and therefore it only remains to consider the integral over the remaining set  $[2x_1,2x_2]\times[0,2x_2]$. 
Let us  denote $M:=\|\nabla\omega(t,\cdot)\|_{L^\infty([0,2x_2]^2)}$ and  write $\omega(t,y)=v(t,y_1)+w(t,y)$, where $v(t,y_1):=\omega(t,y_1,x_2)$ and so 
\[
|w(t,y)|=|\omega(t,y_1,y_2)-\omega(t,y_1,x_2)|\le M|x_2-y_2|.
\]
For $y_1\in[2x_1, 2x_2]$ we have with $z:=y_2-x_2$,
\[
\int_{0}^{2x_2} \frac{y_1(x_2-y_2)}{|x-y|^2 |x-\til y|^2} v(t,y_1)dy_2 = \int_{-x_2}^{x_2} \frac{-y_1z}{[(x_1-y_1)^2 + z^2] [(x_1+y_1)^2+z^2]} v(t,y_1) dz = 0
\]
and we also have
\begin{align*}
\left| \int_{2x_1}^{2x_2} \int_{0}^{2x_2} \frac{y_1(x_2-y_2) w(t,y)}{|x-y|^2 |x-\til y|^2}  dy_2 dy_1 \right| 
  \le CM  \int_{2x_1}^{2x_2} \int_{0}^{x_2} \frac{y_1z^2}{(y_1^2 + z^2)^2} dz dy_1 
 \le CM \int_{0}^{x_2} \frac{z^2}{x_1^2 + z^2} dz \le CMx_2. 
\end{align*}
This gives the second term in the min in \eqref{2.7}.

We note that for the second term in \eqref{2.4}, one can always integrate over $[2x_1,1]\times[0,2x_2]$ because the substitution $z_1:=y_1-x_1$ and $z_2:=y_2+x_2$ yields 
\[
C \int_{x_1}^1 \int_{x_2}^{3x_2}  \frac{z_1 z_2}{(z_1^2+z_2^2)^2} dz_2dz_1 
\le C \int_{x_2}^{3x_2}  \frac{z_2}{x_1^2+z_2^2}  dz_2 
\le C   \log \frac{x_1^2+3x_2^2}{x_1^2+x_2^2} \le C.
\]

\end{proof}

\begin{proof}[Proof of Theorem \ref{T.1.1}]
(i) Given $\alpha$ and $A$, pick a function $\omega_0:(2\bbT)^2\to [-1,1]$ which is  odd in both $x_1$ and $x_2$, non-negative on $[0,1]^2$ and equal to 1 on a subset of $[0,1]^2$ of measure $1-\delta$ (for some $\delta\in(0,\tfrac 1{10})$ to be chosen later), with $\omega_0\in C^{1,\alpha}((2\bbT)^2)\cap C^\infty_{\rm loc}((2\bbT)^2\setminus\{0\})$ and  $\omega_0(s,s)=s^{1+\alpha}$ for $s\in[0,\delta]$. 
For instance, on $B_\delta (0)$ we could have in polar coordinates $\omega_0 (r,\phi)=(r/\sqrt 2)^{1+\alpha} \sin(2\phi)$.

Take any $T\ge T_0:= \tfrac 1A  |\log \delta|$
so that $e^{-AT} \le \delta$, let $X(t)$ solve $X'(t)=u(t,X(t))$ with $X(0)=(e^{-aAT},e^{-aAT})$ for some $a>1$ to be chosen later, and let $T':=\min\{T,T^*\}$, with $T^*$ the exit time of $X$ from the square $[0,e^{-AT}]^2$.  Obviously, 
\beq \lb{2.3}
\omega(t,X(t))= \omega_0(X(0)) = e^{-a(1+\alpha)AT}
\eeq
 for all $t\ge 0$. Let us also assume that
\beq \lb{2.9}
\sup_{t\le T} \|\nabla\omega(t,\cdot)\|_{L^\infty([0,2\exp(-AT)]^2)}\le  e^{AT}
\eeq
because otherwise we are done.  Since $X(t)\in [0,e^{-AT}]^2$ for $t\le T'$, 
\beq\lb{2.6}
x_2 \|\nabla\omega(t,\cdot)\|_{L^\infty([0,2x_2]^2)}\le 1
\eeq
in \eqref{2.7} when  $t\le T'$ and  $x=X(t)$ (and the same estimate applies to \eqref{2.8}).  We then have \eqref{2.1} 
with $|B_j(t,X(t))|\le C$ for $t\le T'$   (recall that $\|\omega(t,\cdot)\|_{L^\infty}=\|\omega_0\|_{L^\infty}=1$).

A crucial observation of \cite{KS} is that $\omega_0\ge 0$ on $[0,1]^2$ and $\omega_0=1$ on a subset of $[0,1]^2$ of measure $1-\delta$ (along with the distribution function of $\omega(t,\cdot)$ being the same for all $t$) guarantees that the integral in \eqref{2.1} is no less than $\tfrac 1C|\log\delta|$ when $\delta<\tfrac 1{10}$ and $x\in[0,\delta]^2$, for some universal $C>0$.  (If instead we had odd $\omega_0:(2\bbT)^2\to [-\eps,\eps]$ equal to $\eps$ on a subset of $[0,1]^2$ of measure $1-\delta$, then this would be $\tfrac \eps C|\log\delta|$, and our proof would be unchanged.)

So if we denote by $k(t)$ the value of the integral in \eqref{2.1} for $x=X(t)$, multiplied by $\tfrac 4\pi$,  then $k(t)\ge \tfrac 1C|\log\delta|$ for $t\le T'$.  Hence we have for $t\le T'$,
\begin{align}
 u_1(t,X(t)) &\le - \left(\frac 1C |\log\delta| -C \right) X_1(t), \lb{2.10}
\\ u_2(t,X(t)) &\ge \left(\frac 1C |\log\delta| -C \right) X_2(t). \lb{2.11}
\end{align}
If we take $\delta<e^{-C^2}$, it follows that $X_1(T')<e^{-AT}$ and  $T'\le (a-1)AT(\tfrac 1C |\log\delta| -C)^{-1}$.  
We will in fact  pick $\delta\le e^{-C(a-1)A-C^2}$  so that also $(a-1)A(\tfrac 1C |\log\delta| -C)^{-1}< 1$.  Hence $T'=T^*< T$ and $X_2(T')=e^{-AT}$. 
In addition, 
\beq\lb{2.12}
\left| \frac d{dt}\left[ \log X_1(t) + \log X_2(t) \right] \right| \le C
\eeq
for $t\le T'$ by Lemma \ref{L.2.1} and \eqref{2.6}.
Therefore,
\beq\lb{2.13}
\log X_1(T') \le \log X_1(0)- \log X_2(T')+\log X_2(0) + CT' \le [-2aA+A+C]T.
\eeq
But this, \eqref{2.3}, and $\omega(t,0,e^{-AT})=0$ give
\[
\log \|\nabla\omega(T',\cdot)\|_{L^\infty([0,\exp(-AT)]^2)} \ge \log \frac{\omega(T',X(T'))}{X_1(T')} \ge 
[a(1-\alpha)A-A-C]T, 
\]
which equals $AT$ if we pick 
\[
a:= \frac{2A+C}{(1-\alpha)A}
\]
and then $\delta$ as above.  The proof of (i) is finished because $T'\le T$.

(ii) Given $A$, pick $\omega_0:(2\bbT)^2\to [-1,1]$ which is  odd in both $x_1$ and $x_2$, non-negative on $[0,1]^2$ and equal to 1 on a subset of $[0,1]^2$ of measure $1-\delta$ (for some  $\delta\in(0,\tfrac 1{10})$ to be chosen later), smooth, and with  $\omega_0(x_1,x_2)=\sin^3 (\pi x_1) \sin (\pi x_2)$ when $\min\{|x_1|,|x_2|\}\le \tfrac \delta 4$. 

Take any $T\ge T_0:= \tfrac 1A  |\log \tfrac \delta 4|$
so that $e^{-AT} \le \tfrac \delta 4$, let $X(t)$ solve $X'(t)=u(t,X(t))$ with $X(0)=(e^{-AT},e^{-(2a-1)AT})$ for some $a>1$ to be chosen later, and let $T':=\min\{T,T^*\}$, with $T^*$ the exit time of $X$ from the square $[0,e^{-AT}]^2$.  Obviously, 
\beq \lb{2.20}
\omega(t,X(t))= \omega_0(X(0)) = \sin^3 (\pi e^{-AT}) \sin (\pi e^{-(2a-1)AT}) \ge e^{-(2a+2)AT}
\eeq
 for all $t\ge 0$. Let us also assume \eqref{2.9} (because otherwise we are done, using $\nabla\omega(t,0)=0$).
 As above, we obtain \eqref{2.1} 
 with $|B_j(t,X(t))|\le C$ for $t\le T'$.
 
We thus again have \eqref{2.10}--\eqref{2.12} for $t\le T'$, as well as  $X_2(T')=e^{-AT}>X_1(T')$ and $T'=T^*< T$, provided we pick $\delta< e^{-C(2a-2)A-C^2}$.  So \eqref{2.13} follows as well and then by \eqref{2.20},
\beq\lb{2.21}
\log \sup_{s\in[0,\exp([-2aA+A+C]T)]} \omega_{x_1}(T',s,e^{-AT}) \ge \log \frac{\omega(T',X(T'))}{X_1(T')} \ge 
-[3A+C]T. 
\eeq
The result will follow if we can show that $\omega_{x_1}(T',0,e^{-AT})=0$, because this and \eqref{2.21} imply
\[
\log \|D^2 \omega(T',\cdot)\|_{L^\infty([0,\exp(-AT)]^2)} \ge [2aA-4A-2C]T,
\] 
provided $a\ge 1+\tfrac C{2A}$ (so that $-2aA+A+C\le -A$).  Then we only need to pick
\[
a:=\frac{5A+2C}{2A}.
\]

Let $v(t,x):=\omega_{x_1}(t,x)$.  Then
\[
v_t+u_1 v_{x_1} + u_2 v_{x_2} + (u_1)_{x_1} v + (u_2)_{x_1} \omega_{x_2} = 0. 
\]
We have $u_1(t,0,x_2)=0=\omega(t,0,x_2)$ by symmetry, so also $\omega_{x_2}(t,0,x_2)=0$.
This shows that if we denote $w(t,s):=v(t,0,s)$ for $s\in 2\bbT$, then for $(t,s)\in\bbR^+\times 2\bbT$,
\[
w_t+ u_2(t,0,s) w_{s} + (u_1)_{x_1}(t,0,s) w = 0. 
\]
Since $w(0,s)=(\omega_0)_{x_1}(0,s)= 0$ and  $u$ is smooth, it follows that $w\equiv 0$.  Thus we indeed obtain $\omega_{x_1}(T',0,e^{-AT})=0$, and the proof of (ii) is finished.
\end{proof}




\end{document}